\documentclass[10pt]{article}

\usepackage{bm}
\usepackage[top=0.9in,bottom=1.1in,left=1.05in,right=1.05in]{geometry}

\usepackage{textgreek,mathtools}
\usepackage{amsmath}
\usepackage{graphicx,subfig}
\usepackage[colorinlistoftodos]{todonotes}
\usepackage[colorlinks=true, allcolors=blue, urlcolor=black]{hyperref}
\usepackage{amsmath}
\usepackage{upgreek}
\usepackage{amssymb}
\usepackage{amsfonts}
\usepackage{amsthm}
\usepackage{mathrsfs}
\usepackage{fontenc}
\usepackage{empheq}
\usepackage{enumerate}
\usepackage{comment}
\usepackage{appendix}
\usepackage{bbm}
\usepackage{algorithm,algpseudocode}
\mathtoolsset{showonlyrefs}

\theoremstyle{plain}
\newtheorem{thm}{Theorem}[section]

\theoremstyle{definition}

\theoremstyle{remark}

\newcommand{\email}[1]{\protect\href{mailto:#1}{#1}}

\title{Empowering Optimal Control with Machine Learning:\\ A Perspective from Model Predictive Control}
\date{\today}
\author{
Weinan E\thanks{AI for Science Institute, Beijing; 
         Center for Machine Learning Research
         and School of Mathematical Sciences, 
        Peking University, Beijing, China (\email{ weinan@math.pku.edu.cn}).}
\and Jiequn Han\thanks{Center for Computational Mathematics, 
        Flatiron Institute, New York, 10010, USA (\email{jiequnhan@gmail.com})}
\and Jihao Long\thanks{Program of Applied and Computational Mathematics,
          Princeton University,
          Princeton, 08544, USA (\email{ jihaol@princeton.edu})}
}

\begin{document}
\maketitle

\begin{abstract}
Solving complex optimal control problems have confronted computational challenges for a long time. Recent advances in machine learning have provided us with new opportunities to address these challenges.
 This paper takes model predictive control, a popular optimal control method, as the primary example to survey recent progress that leverages machine learning techniques to empower optimal control solvers.
 We also discuss some of the main challenges encountered when applying machine learning to develop more robust optimal control algorithms.
\end{abstract}

\section{Introduction}

Control algorithms\,\cite{kirk2004optimal} are  widely used in engineering and industry. 
Its objective is to find an optimal control for a dynamical system that minimizes a given loss function.
The so-called open loop control aims to address this problem for a pair of specific initial and terminal conditions.
 Mathematically speaking, this requires solving a two-point boundary value problem connecting
the initial and terminal states.  If for some reason the state of the system deviates from the designed path, one has to 
redo all the calculations to find the new optimal control.
Closed loop control, on the other hand, aims to find the optimal policy function: the optimal action as a function of the
state. This optimal policy function can be obtained from the solution of the associated Bellman equation.
Compared to the open loop control, closed loop control is more useful in practice thanks to its flexibility when the initial state varies.
However, since the time of Bellman, it has been realized that solving high-dimensional closed loop control problems is a difficult task \cite{bellman1957}.
The popularly used terminology ``curse of dimensionality''  was originally coined in order to highlight these difficulties \cite{bellman1957}.
For that reason, theoretical work on optimal control algorithms has been limited to either the open loop control
or low-dimensional closed loop control problems, and one often has to resort to uncontrolled approximations, such as assuming the value or policy functions admitting specific low-dimensional structures, to address practical needs \cite{powell2007approximate}.

 The success of machine learning in computer vision and other AI tasks has raised new hope.
 In essence, what deep learning does in image recognition, image generation and AlphaGo \cite{silver2016mastering} is to solve some underlying 
 mathematical problems such as  approximating functions, approximating probability distributions and the Bellman equation respectively
 \cite{goodfellow2016deep,e2021mlmodeling}. These mathematical problems have been the subject of computational mathematics and numerical analysis for a very long time, except that the problems treated by deep learning have much higher dimensionality than the ones treated in computational mathematics. 
 It is then natural to ask whether similar techniques can be used to address high-dimensional problems in scientific computing
 such as optimal control problems.
 The very first attempt in this direction was made in \cite{han2016deep}, which  demonstrated that deep learning-based algorithms can be used to solve, fairly routinely, 
 stochastic control problems whose dimensionality can be  as large as 50.
 Since then the field of high-dimensional control problems has flourished. 
 This paper aims to review recent progress on machine learning enhanced optimal control algorithms.
 We will also  discuss  the new challenges that arise in this exciting new area.

We start with a brief review of the two most popular  control algorithms used in practice, the proportional-integral-derivative (PID) controller \,\cite{aastrom1995pid} and model predictive control (MPC)\,\cite{richalet1978model,cutler1980dynamic}.
A PID controller continuously calculates the error between the reference point and the measured output of the system. The control then corrects the error through a linear combination of the error itself, the time integral of the error, and the time derivative of the error with some predetermined coefficients. Although the PID controller may not provide an optimal solution, it has been the most popular choice in industrial applications\,\cite{samad2020industry} due to its high feasibility and reasonably satisfactory performance. However, its performance  severely degrades if the system has strong nonlinearity, time-inhomogeneity, or multiple inputs/outputs. It  also has difficulties dealing with constraints. 

MPC on the other hand, is developed to handle the situation when the time horizon for the control problem is large.
Problems with a long time horizon might suffer from an additional curse, much like the ``curse of memory'' found in the context of learning
dynamical systems \cite{li2020curse,li2022approximation}, even in low dimension.
 MPC solves a sequence of (auxiliary) online optimal control problems with short time horizons. %
At each time step $t$, the MPC solver receives a measured state and solves an optimal control problem whose state trajectory starts from the measured state at $t$ and %
is defined on a shorter time interval.
The controls obtained  are then applied to the system until new measurements are 
received. 
Unlike the PID controller, MPC explicitly utilizes the model of the dynamical system and hence is able to handle complex systems with many inputs, outputs and constraints. For this reason, it is widely used in process industry\,\cite{qin2003survey}, power electronics\,\cite{vazquez2014model,kouro2015model,karamanakos2020model}, robotics\,\cite{duchaine2007computationally,nubert2020safe,kleff2021high}, unmanned aerial vehicles\,\cite{kang2009linear,kamel2017model,altan2020model}, etc.. 
However, the computational cost has been an obstacle in the application of MPC in high-frequency control problems.
In addition, as we will see below, the performance of MPC highly depends on the quality of the reference trajectory, which is also time-consuming to generate. These challenges are aggravated in high-dimensional nonlinear problems. 

Overall, classical control algorithms like PID and MPC face computational
 challenges for high-dimensional nonlinear problems where machine learning techniques 
demonstrate a clear advantage. This advantage motivates us to apply machine learning techniques to high-dimensional control problems. We will review recent progress and discuss new challenges that arise in this exciting new area.
Our main focus is on MPC since the issues encountered there are quite representative of the issues we are interested in.
 We do mention other control algorithms when they are of relevance.

\section{A Short Review of Model Predictive Control}
We work with the discrete-time optimal control problem:
\begin{align}\label{OCP_Problem}
    &\min_{\{\bm x_k,\bm u_k\}_{k=0}^T} \sum_{k=0}^{T-1} L_k(\bm x_k,\bm u_k) +  M(\bm x_T) \\
    \text{s.t. } &\begin{cases}
    \bm x_{k+1} = \bm f_k(\bm x_k,\bm u_k),\, 0 \le k \le T-1, \bm x_0 = \bm x\\
    \bm g_k(\bm x_k,\bm u_k) = 0,\; \bm h_k(\bm x_k,\bm u_k) \ge 0,\; 0 \le k \le T,
    \end{cases}\notag
\end{align}
in which $L_k$ is the running cost, $M$ is the terminal cost, and $\bm g_k, \bm h_k$ denote state-action constraints. 
Throughout this paper, a regular character indicates a scalar, a boldface character indicates a vector or matrix, and $\ge$ means component-wise order.
We call \eqref{OCP_Problem} global optimal control problem in contrast to the local optimal control problem used in MPC.
Since continuous-time optimal control problems can be discretized through various discretization schemes, all the discussions below also apply to continuous-time optimal control problems. 

\begin{figure}[ht]
  \center
  \includegraphics[width=8.4cm]{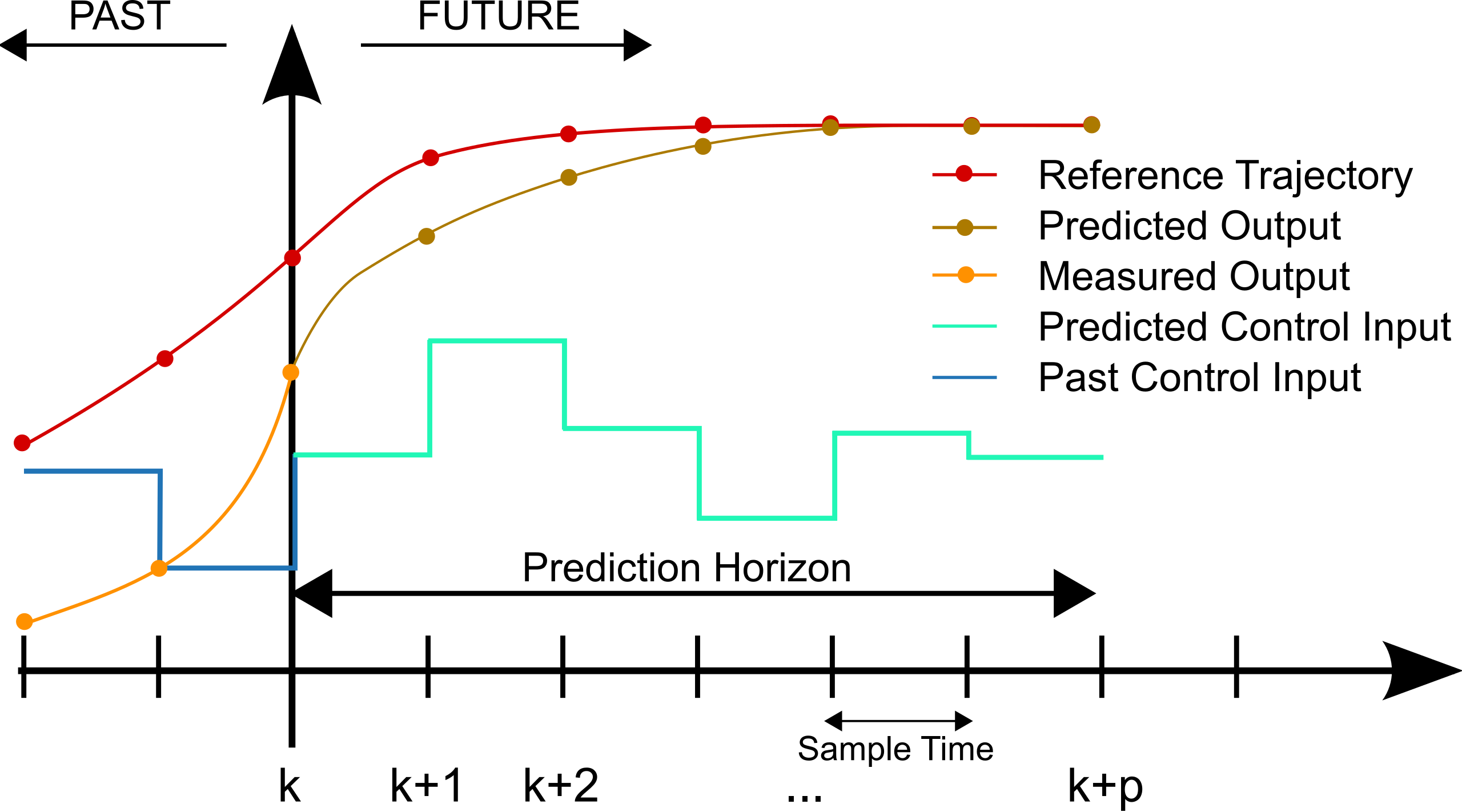}
  \caption{A schematic plot of an MPC solver, adapted from~\cite{MPC_wiki}.}
  \label{fig:mpc_scheme}
\end{figure}

We are interested in solving \eqref{OCP_Problem} for relatively large  time horizon $T$.  As mentioned above, this is in general a challenging
problem. Moreover, 
due to the accumulation of the numerical error and model misspecification error, the control obtained through solving \eqref{OCP_Problem} may not perform well in practice. The idea of  model predictive control (MPC) is to approximately solve the global optimal control problem \eqref{OCP_Problem} by solving a sequence of online local optimal control problems. 
At each sample time $t$, the current state of the system  $\bm x$ is measured.  With that one  considers the following local optimal control problem:

\begin{align}\label{MPC_Problem}
    &\min_{\{\bm x_k, \bm u_k\}_{k=t}^{t_p}} \sum_{k=t}^{t_p-1} L_k(\bm x_k,\bm u_k) +  \tilde{M}(t_p,\bm x_{t_p}) \\
    \text{s.t. } &\begin{cases}
    \bm x_{k+1} = \bm f_k(\bm x_k,\bm u_k),\, t \le k \le t_p-1, \bm x_t = \bm x\\
    \bm g_k(\bm x_k,\bm u_k) = 0,\; \bm h_k(\bm x_k,\bm u_k) \ge 0,\; t \le k \le t _p,
    \end{cases}\notag
\end{align}
where $t_p = \min\{t+p,T\}$,  
the prediction horizon $p$  is much smaller than $T$, and $\tilde{M}$ is a heuristic approximation of the value function $V$:
\begin{align*}
    V(t',\bm x') &= \min_{\{\bm x_k,\bm u_k\}_{k= t'}^{T}} \sum_{k=t'}^{T-1} L_k(\bm x_k,\bm u_k) + M(\bm x_T) \\
    \text{s.t. } &\begin{cases}
    \bm x_{k+1} = \bm f_k(\bm x_k,\bm u_k),\, t' \le k \le T-1, \bm x_{t'} = \bm x'\\
    \bm g_k(\bm x_k,\bm u_k) = 0,\,\bm h_k(\bm x_k,\bm u_k) \ge 0,\, t' \le k \le T.
    \end{cases}
\end{align*}
After the local problem \eqref{MPC_Problem} is solved, the controls obtained are applied to the system until the next sampling.
At the next sampling time, the current state is updated  and the same procedure is repeated. %
Since the prediction horizon $p$ is much smaller than the terminal time $T$, solving the local  problem \eqref{MPC_Problem} is much easier than solving the global  problem \eqref{OCP_Problem}. 
The online interaction between the MPC solver and the system significantly reduces the accumulation of the numerical error and model misspecification error.  The difficulty of this MPC formulation lies in obtaining an accurate  approximation of the value function.
A common remedy in practice is to use a reference trajectory. Namely, let $\{\bm y_k\}_{k=0}^{T-1}$ be the reference trajectory and set
\begin{equation*}
    L_k(\bm x, \bm u) = \|\bm y_k - \bm r_k(\bm x, \bm u)\|^2, \,\tilde{M} =0.
\end{equation*}
The reference trajectory is produced by human expertise or other planning algorithms\,\cite{lavalle2006planning} to ensure that the solution of MPC has good global performance, i.e. the total cost of the MPC solution in the global optimal control problem \eqref{OCP_Problem} is small. Figure \ref{fig:mpc_scheme} gives a schematic plot on how a typical MPC solver works.

In traditional MPC, \eqref{MPC_Problem} is formulated as a linear quadratic control problem and the closed-form solution is obtained by solving the Riccati differential equation \cite{anderson2007optimal}. For nonlinear MPC, various algorithms including direct shooting method\,\cite{bock1984multiple}, direct collocation method\,\cite{logsdon1989accurate}, indirect shooting method\,\cite{keller1976numerical}, indirect collocation method\,\cite{russell1972collocation} and dynamic differential programming \,\cite{mayne1966second,jacobson1970differential} can be used to solve \eqref{MPC_Problem};  see\, \cite{betts1998survey,rao2009survey} for detailed discussions.

MPC is a closed loop controller since it implicitly defines a policy function between the state $\bm{x}_t$ and control $\bm{u}_t$ at sample time $t$ through the local optimal control problem \eqref{MPC_Problem}. However, unlike usual closed loop controllers, the computation of problem \eqref{MPC_Problem} is online and one needs to redo the computation if the initial state changes. Moreover, if the initial state changes, the reference trajectory should be updated as well. In this sense, MPC using the reference trajectory is more like an open loop controller, which introduces several challenges when dealing with practical problems.
One challenge is the computational cost. Since MPC is an online algorithm, the computational time for a single optimization problem \eqref{MPC_Problem} should not be larger than the sampling interval for updating the states. Since the prediction horizon $p$ significantly impacts the computational time, this requirement limits the scale $p$.
In some areas such as power electronics and signal processing\,\cite{geyer2008model}, the interval between the sample time is so small that the prediction horizon is usually set to one. On the other hand, the prediction horizon $p$ is also an essential factor for the performance of MPC.
Ideally $p$ should be large enough to fully take into account the influence of the current control on the whole dynamics.
Hence, reducing the computational time for relatively large prediction horizons is essential for many applications, especially for high-dimensional nonlinear problems.

Besides the computational cost, another challenge is providing a good reference trajectory efficiently. Even with a relatively large prediction horizon, the performance of MPC still hinges on the reference trajectory. The planning algorithms that produce reference trajectories must be re-run if the initial states or other parameters of the control problems change. How to efficiently generate reference trajectories of high-quality or design other methods to provide global guidance still remains an open question.

\section{Empowering MPC with Machine Learning}

\begin{figure}[ht]
    \center
  \includegraphics[width=8.4cm]{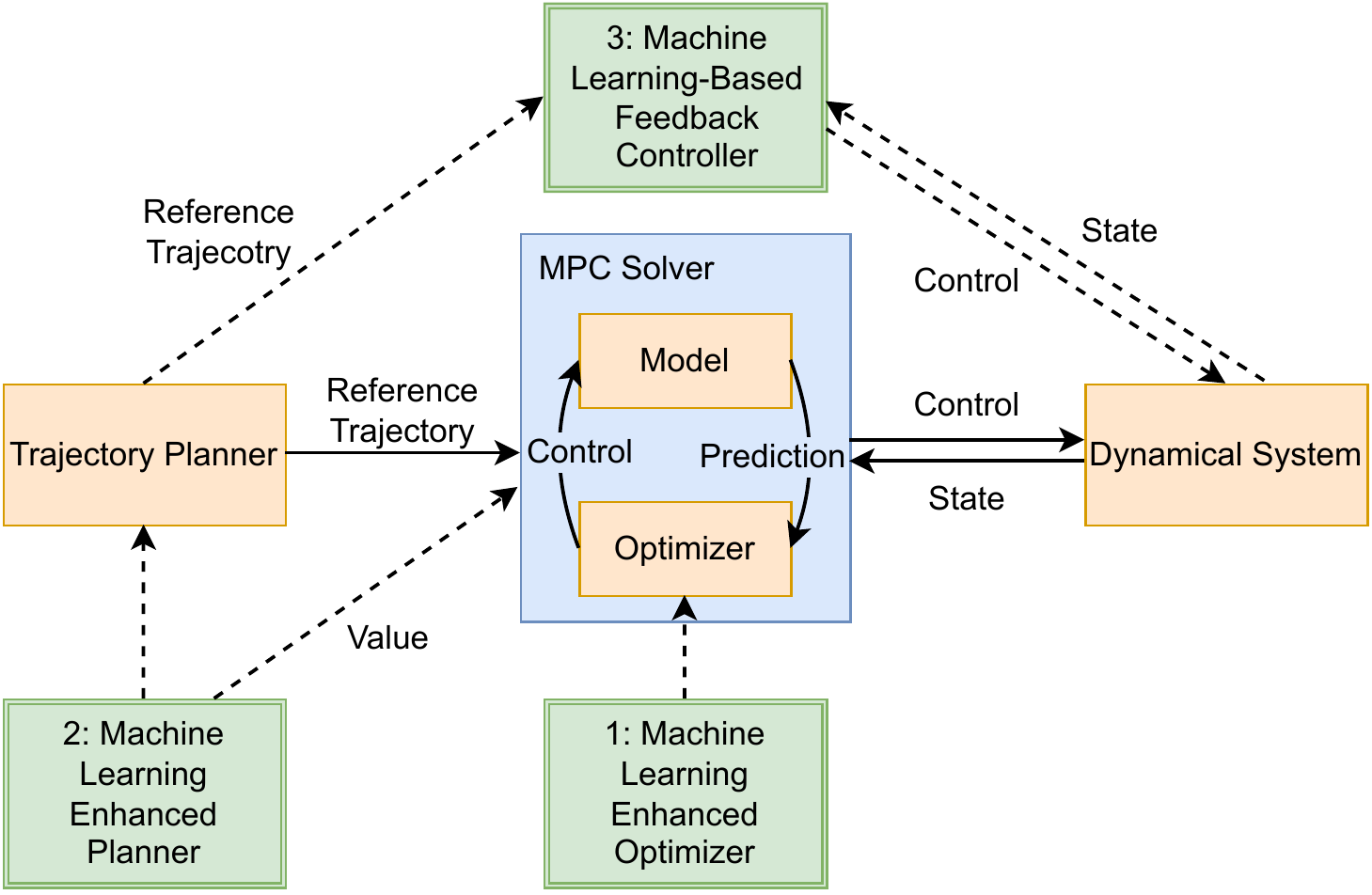}
  \caption{Three ways to empower MPC with machine learning techniques: (1) machine learning enhanced optimizer, (2) machine learning enhanced planner, (3) machine learning-based feedback controller.}
  \label{fig:MPC_ML}
\end{figure}

We now discuss some recent progress that leverages machine learning techniques to tackle the aforementioned computational challenges. Figures \ref{fig:MPC_ML} summarize different aspects to empower MPC with machine learning techniques. As mentioned above, an MPC solver takes the reference trajectory generated by a trajectory planner and measured current state in the dynamical system as the inputs, uses the model to predict the trajectory in the prediction horizon and computes the local optimal trajectory for \eqref{MPC_Problem}, and finally outputs the control of the first step for the dynamical system. Machine learning techniques can improve
 the whole process in at least three ways. First, machine learning techniques can accelerate the optimizer in the MPC solver. Second, machine learning techniques can enhance the trajectory planner or provide an accurate approximation to the value function, both of which can improve the global performance of MPC. Finally, we can use machine learning to approximate the feedback control, in the form of either the optimal closed loop control or the MPC, to replace the MPC solver. We detail the three ways in the following discussions.

\textbf{Using machine learning to accelerate solving the MPC problem.} 
\,\cite{nakamura2021adaptive} proposes the neural-network warm start technique to solve such optimal control problems. In \cite{nakamura2021adaptive}, training data is generated by solving the optimal control problem \eqref{MPC_Problem} with various initial conditions and then used to train a neural network to approximate the value function offline. 
The trained functions can then generate a path that is usually reasonably close to the optimal path and hence serves as a good initial guess for the optimal control solver in online MPC. Note that although the work of \cite{nakamura2021adaptive} is based on the indirect collocation method, this is a general idea applicable for almost all nonlinear optimal control solvers. Besides providing the whole path for the initial guess, neural networks can also provide the crucial parameters in the optimal control problem and accelerate the computation. 
\cite{zang2022machine} presents a neural network-based algorithm for  learning the optimal terminal time for a landing problem. With the optimal terminal time accurately estimated by the neural network, the optimal control problem becomes much easier to solve: the 
 success rate increases and computational time decreases. 
 In addition, neural networks can also be used to tune the critical parameters in MPC to improve the performance, see, e.g., \cite{dragivcevic2018weighting}.

\textbf{Using machine learning to enhance the global planner.}\,
As mentioned above, MPC only produces locally optimal policies, and hence a planning algorithm that provides global guidance is vital for the performance of MPC. The reference trajectories provided by traditional planning algorithms are hard to reuse if the initial states or other parameters change. \cite{zhong2013value} proposes an algorithm to learn the value function offline and replace the heuristic approximation $\tilde{M}$ with the approximated value function in the optimal control problem \eqref{MPC_Problem}. The dynamic programming principle ensures the global performance of the resulting algorithms.
Compared to other planning algorithms, the value function can be easily reused when the initial conditions change. In \cite{zhong2013value}, nearest neighbor, locally weighted projection regression, and polynomial mixture of normalized Gaussians are used as function approximators. 
Subsequent  works \cite{deits2019lvis,landry2021seagul} extended these to neural networks.

Since many planning algorithms are related to solving the global optimal control problems~\eqref{OCP_Problem}, the techniques used to accelerate solving the MPC problem~\eqref{MPC_Problem} can also be applied to~\eqref{OCP_Problem}, see \cite[Section~4]{zang2022machine} for an example. 

\textbf{Using machine learning to approximate the feedback control.}\,
Since the original work in \cite{han2016deep}, there has been a flurry of activity on developing deep learning-based algorithms
for  approximating the feedback control. The training data is either generated by MPC\, \cite{hertneck2018learning,stevvsic2018sample,mohamed2019neural,nubert2020safe,wu2021composing} or the global optimal control method\,\cite{han2016deep,zhang2019safe,nakamura2020qrnet,nakamura2021adaptive,nakamura2021neural}. There are several advantages to training a machine learning model to approximate the MPC or optimal policies offline. The deployment is straightforward, the online computational time is extremely small, and the same model can be used with different initial conditions or other model parameters. Nevertheless, as pointed out in \,\cite{nakamura2021adaptive,zang2022machine}, neural network feedback controllers can fail to stabilize a system, even when they are trained to a high degree of accuracy.  This phenomenon may be related to the adversary attack problem of neural networks\,\cite{szegedy2013intriguing}. Some strategies are proposed to address this issue.  \cite{nakamura2021neural} discussed adjusting the structure of the neural networks to imitate  linear quadratic control around the stationary point.  \cite{wu2021composing}  suggested that if the state is near
 the stationary point, then linear quadratic control can be used; otherwise, one can use MPC periodically and the neural-network controller at other times. In \cite{zhang2019safe}, a validation criterion is defined, and a neural network controller is only used if the criterion is satisfied; otherwise, a backup controller such as MPC is used. Although these strategies can improve the stability of the neural network controller, there is still room
 for improvement.

\section{Discussions}
As discussed above, offline trained machine learning models for approximating the value or policy functions can help accelerate solving the MPC problem, enhance the global planner, and even directly be used to control the dynamics.
Existing results have shown the potential of machine learning empowered optimal control solvers. 
However, this field is still developing and faces many challenges. Here we list three challenges that, in our view, are crucial.
\begin{enumerate}
    \item \textbf{How to generate the training data?} The training data  usually consists of optimal or near-optimal paths of the control problems with various initial states and other parameters. To generate the training data, one needs to solve a class of open loop control problems. 
    Ideally one would like to generate the optimal training dataset that is  both small and representative enough for the practical situations of interest.
    Despite much progress, data generation is still hard for practical problems such as unmanned aerial vehicles moving in a complex environment. The neural-network warm start technique proposed  in \,\cite{nakamura2021adaptive} and other related techniques demonstrate the possibility of using machine learning to help with the data generation process.  But there is a lot of room for improvement in
    this direction.
      \item \textbf{How to efficiently train the machine learning model?}  Training is always crucial for neural network models.
      In the current task, since data generation is part of the model training process, we have the additional issue of
       interactively generating the data and training the model. 
     Some initial progress has been made in \cite{nakamura2021adaptive}.
     We believe that the exploration-labeling-training
     (ELT)  algorithms \cite{zhang2018reinforced,e2021mlmodeling} can guide the interaction between generating the data and training the model, which proceeds iteratively with the following steps: (1) exploring the state space and examining which states need to be labeled; (2) computing the control problem for labeled states and placing them into the training data; (3) training the machine learning model. We look forward to a systematic study on ELT algorithms for optimal control problems in the near future.
    \item \textbf{How to use the approximate value or policy function?}  Measures for the performance of optimal control algorithms include computational time, optimality, stability, and robustness against model misspecification. While machine learning empowered algorithms often outperform  traditional algorithms in the first two aspects, existing work has not  considered the last two aspects much. How to design algorithms that not only inherit the advantages of computational time and accuracy but also guarantee  stability and robustness is a very important
    issue that should be addressed in the near future.
\end{enumerate}

\bibliographystyle{plain}
\bibliography{Reference}

\end{document}